\documentclass[11pt]{article}
\usepackage{amsfonts}

\newtheorem{thm}{Theorem}

\usepackage{epsfig}

\begin{document}

 \def\d{{\delta}}
 \def\ci{{\circ}}
 \def\e{{\epsilon}}
 \def\l{{\lambda}}
 \def\L{{\Lambda}}
 \def\m{{\mu}}
 \def\n{{\nu}}
 \def\o{{\omega}}
 \def\s{{\sigma}}
 \def\v{{\varphi}}
 \def\a{{\alpha}}
 \def\b{{\beta}}
 \def\p{{\partial}}
 \def\r{{\rho}}
 \def\ra{{\rightarrow}}
 \def\lra{{\longrightarrow}}
 \def\g{{\gamma}}
 \def\D{{\Delta}}
 \def\La{{\Leftarrow}}
 \def\Ra{{\Rightarrow}}
 \def\x{{\xi}}
 \def\c{{\mathbb C}}
 \def\z{{\mathbb Z}}
 \def\2{{\mathbb Z_2}}
 \def\q{{\mathbb Q}}
 \def\t{{\tau}}
 \def\u{{\upsilon}}
 \def\th{{\theta}}
 \def\la{{\leftarrow}}
 \def\lla{{\longleftarrow}}
 \def\da{{\downarrow}}
 \def\ua{{\uparrow}}
 \def\nwa{{\nwtarrow}}
 \def\swa{{\swarrow}}
 \def\nea{{\netarrow}}
 \def\sea{{\searrow}}
 \def\hla{{\hookleftarrow}}
 \def\hra{{\hookrightarrow}}
 \def\sl{{SL(2,\mathbb C)}}
 \def\ps{{PSL(2,\mathbb C)}}
 \def\qed{{\hfill$\diamondsuit$}}
 \def\pf{{\noindent{\bf Proof.\hspace{2mm}}}}
 \def\ni{{\noindent}}
 \def\sm{{{\mbox{\tiny M}}}}
 \def\sc{{{\mbox{\tiny C}}}}
 \def\ke{{\mbox{ker}(H_1(\p M;\2)\ra H_1(M;\2))}}

\begin{center}
{\bf Heegaard Splittings and Virtually Haken Dehn Filling}
\end{center}

\begin{center}
J. Masters\footnote{{Supported by NSF Postdoctoral Fellowship}},
 W. Menasco and X. Zhang\footnote{{Partially supported by NSF
grant DMS 0204428}}
\end{center}

\noindent{\bf Abstract}. We use Heegaard splittings to give some
examples of  virtually Haken 3-manifolds.

\vspace{10mm}
 A compact connected  3-manifold is said to be
virtually Haken if it has a finite sheeted covering space which is Haken.
 The virtual Haken conjecture states
 that every compact, connected, orientable,
 irreducible 3-manifold  with infinite fundamental group
 is virtually Haken.
 Since virtually Haken 3-manifolds and Haken 3-manifolds
 possess similar properties, such as
 geometric decompositions and, in the closed case,
 topological rigidity,
 the resolution of this conjecture would provide solutions to
 several fundamental problems about compact 3-manifolds with
 infinite fundamental groups.

  Some  recent  results in attacking the conjecture can be
found in [CL] [BZ] [M] [DT]. A summary of earlier results
 can be found in [K, Problem 3.2].
 For connections between the virtual Haken conjecture,
 Heegaard splittings, and the Property $\tau$ conjecture,
 see [L].

 Motivated by the work of Casson and Gordon ([CG]), we shall show that
 lifted Heegaard surfaces can often be compressed
 to become essential.
 Our techniques can be used to produce many families of
 non-Haken but virtually Haken 3-manifolds, a few of which are
 given here to illustrate the method.
 A more general result will be proved in a forthcoming paper.

 The first named author wishes to thank Cameron Gordon
 for many useful conversations, and the University of Texas
 at Austin for its hospitality.

We proceed to give the examples.
 Let $K_{2n+1}$ be the twist knot in $S^3$ as shown in Figure 1.
Let $M_n$ be the exterior of $K_{2n+1}$, with standard
meridian-longitude framing on $\p M_n$. Recall that a  connected,
compact, orientable 3-manifold whose boundary is a torus is  called
{\it small} if every closed, orientable, embedded,
 incompressible surface is parallel to the boundary, and called {\it
large} otherwise.

\vspace{3mm}
\begin{figure}[!ht]
{\epsfxsize=3in \centerline{\epsfbox{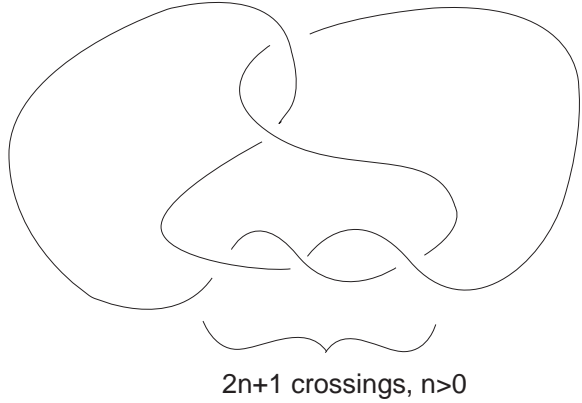}}\hspace{10mm}}
\caption{The twisted  knot $K_{2n+1}$}
\end{figure}

\begin{thm}\label{3fold}
 The $3$-fold cyclic cover of $M_n$
is large for every  $n>0$. Every Dehn filling of $M_n$  with slope
$3p/q$, $(3p,q)=1$, $|p|>1$, yields a virtually Haken 3-manifold.
\end{thm}

 Note that by [HT], $M_n$ is hyperbolic, small,
 and has exactly three boundary slopes, for every $n>0$.
 It follows (combining with [CGLS, Theorem 2.0.3]) that all but
 exactly three  Dehn fillings of $M_n$
 give irreducible non-Haken 3-manifolds. Also note that
each $K_{2n+1}$, $n>0$, is a non-fiberd knot with a genus one
Seifert surface,
  and thus  by [CL] it was
known that  every $m$-fold cyclic cover of $M_n$, $m\geq 4$, is
large and every Dehn filling of $M_n$ with slope $p/q$, $(p,q)=1$,
$|p|\geq 8$, is virtually Haken.

\pf Let $\tilde M_n$ be the 3-fold cyclic cover of $M_n$ with
induced meridian-longitude framing on $\p \tilde M_n$. We shall
show that $\tilde M_n$  contains a connected, essential (i.e.
orientable, incompressible, non-boundary-parallel) genus two closed
surface which has an essential simple closed curve isotopic to a
longitude curve of the cover. It follows from [CGLS, Theorem 2.4.3] that the
surface remains incompressible in every Dehn filling of $\tilde
M_n$ with slope $p/q$, $(p,q)=1$, $|p|>1$.
 As every Dehn filling of $M_n$  with slope
$3p/q$, $(3p,q)=1$, $|p|>1$, is free covered by Dehn filling of
$\tilde M_n$ with slope $p/q$, $(p,q)=1$, $|p|>1$, the second
conclusion of the theorem will follow.

\vspace{3mm}
\begin{figure}[!ht]
{\epsfxsize=5in \centerline{\epsfbox{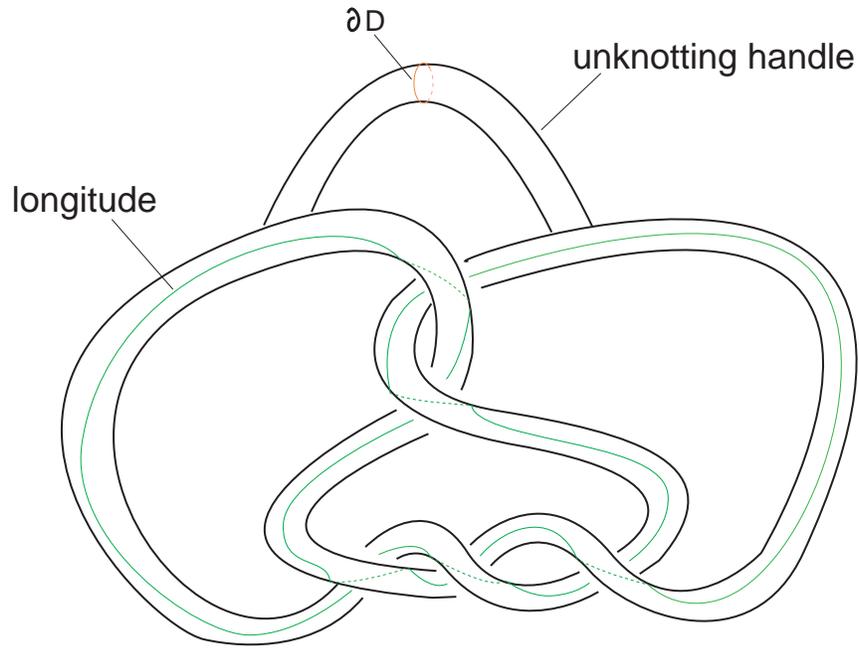}}\hspace{10mm}}
\caption{An unknotting tunnel, its co-core $\p D$ and a standard
longitude of $K$}
\end{figure}

To make the illustration simple, we first prove the theorem with all
details in case $n=1$, i.e. for the $5_2$ knot $K=K_3$.
 The knot $K$ is tunnel number one, and Figure 2 shows an unknotting tunnel.
 Also pictured in Figure 2 is a longitude $\l$ of $K$.
 Let $N$ be a regular neighborhood of $K$ in $S^3$,
$M=M_1=\overline{S^3-N}$, $B$ a regular neighborhood of the
unknotting tunnel in $M$,  and $H=\overline{M-B}$. Then $H$ is a
handle body of genus two.  Let $D$ be a meridian disk of the 1-handle $B$ whose
boundary is shown in Figure 2. We deform the handle body
$H'=N\cup B$ by an isotopy in $S^3$ so that its exterior $H$ can
be recognized as a standard handle body in $S^3$ and at the same
time we trace the corresponding deformation of  $\p D$ and $\l$
under the isotopy. The process is shown through Figures 3-6.

\begin{figure}[!ht]
{\epsfxsize=5in \centerline{\epsfbox{f3.ai}}\hspace{10mm}}
\caption{The deformation of $H'$, $\p D$ and $\l$ (part a)}
\end{figure}

\vspace{3mm}
\begin{figure}[!ht]
{\epsfxsize=5in \centerline{\epsfbox{f4.ai}}\hspace{10mm}}
\caption{The deformation of $H'$, $\p D$ and $\l$ (part b) }
\end{figure}

\vspace{3mm}
\begin{figure}[!ht]
{\epsfxsize=5in \centerline{\epsfbox{f5.ai}}\hspace{10mm}}
\caption{The deformation of $H'$, $\p D$ and $\l$ (part c)}
\end{figure}

\vspace{3mm}
\begin{figure}[!ht]
{\epsfxsize=5in \centerline{\epsfbox{f6.ai}}\hspace{10mm}}
\caption{The deformation of $H'$, $\p D$ and $\l$ (part d)}
\end{figure}

A \textit{meridian
disk system } of a handlebody of genus $g$ is a set of $g$
properly embedded mutually disjoint disks in the handle body such
that cutting the handlebody along these disks results in a
3-ball.
Let $\{X, Y\}$ be a meridian disk system of $H$
  whose boundary are shown in Figure 6. Following $\p D$ in
the given orientation, we get a geometric presentation of the
fundamental group $\pi_1(M)$ of $M$: $$\pi_1(M)=<x,y;
x^{-1}y^{-1}x^{-1}yxyx y^{-1}x^{-1}y^{-1}xyxy>,$$ where $x$ is
chosen such that it has a representative curve which is a simple
closed curve in $\p H$ which is disjoint from $\p Y$ and
intersects $\p X$ exactly once and $y$ is also chosen similarly.
(We shall call such generators {\it dual to} the disk system.)
Also we can read off the longitude in terms of these two
generators: $$\l=yxyx^{-1}y^{-1}x^{-1}y^{-2}x^{-1}y^{-1}
x^{-1}yxyx^2.$$

 Cutting  $H$ along $X$ and $Y$, we get a 3-ball. Figure
7 shows the boundary 2-sphere of the 3-ball, which records $X^+$,
$X^-$, $Y^+$, $Y^-$ and $\p D$. Figure 8 shows $H$ in a standard
position, and $\p D$ in $\p H$.

The exterior of $H$ in $M$ is a compression body which we denote
by $C$. Topologically, $C$ is $\p M\times [0,1]$ with a $1$-handle
attached on $\p M\times \{1\}$. It has two boundary components:
one is $\p M=\p M\times \{0\}$ and the other is the genus two
surface $\p H$. We have that $H\cup _{\p H} C$ is a Heegaard splitting of $M$.

 \vspace{3mm}
\begin{figure}[!ht]
{\epsfxsize=3in \centerline{\epsfbox{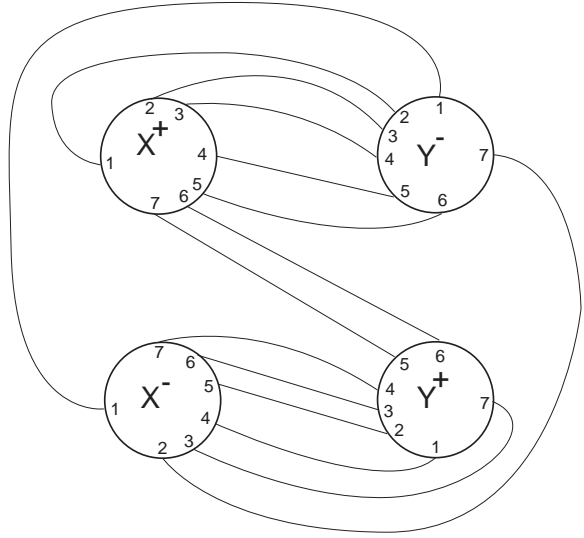}}\hspace{10mm}}
\caption{$\p D$ on the sphere $\p(\overline{H-\{X\times I\cup
Y\times I\}})$ }
\end{figure}

\vspace{3mm}
\begin{figure}[!ht]
{\epsfxsize=5in \centerline{\epsfbox{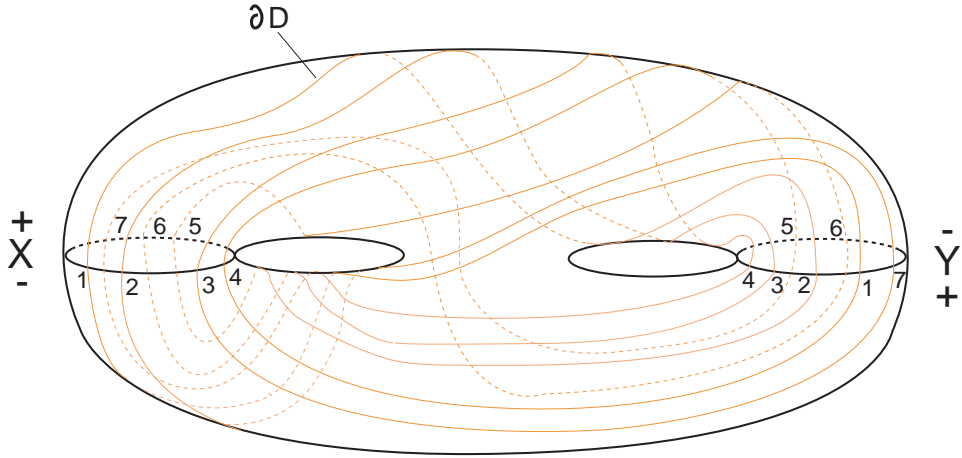}}\hspace{10mm}}
\caption{$H$ and $\p D$ in standard position }
\end{figure}

  Let $\tilde M=\tilde M_1$ be the 3-fold cyclic cover of $M=M_1$.
Note that each of $x$ and $y$ is a generator of $H_1(M;\z)=\z$.
Let $\tilde M$ have the induced Heegaard splitting from that of
$M$. We can easily give the Heegaard diagram of $\tilde M$,  as
shown in Figure 9. The genus four handle body $\tilde H$ in Figure
9 is the corresponding cover of $H$. The corresponding cover
$\tilde C$ of $C$ is a compression body obtained by attaching
three 1-handles to $\p \tilde M\times [0,1]$ on the side $\p\tilde
M\times \{1\}$. The disk $X$ lifts to three disks  $X_1, X_2,
X_3$; and the disk $Y$ lifts to three disks $Y_1, Y_2, Y_3$, as
shown in Figure 9. Pick the meridian disk $X_4$ of $\tilde H$ as
shown in Figure 9. Then \{$X_1$, $X_2$, $X_3$, $X_4$\} forms a
disk system of $\tilde H$. The disk $D$ lifts to three disks
\{$W_1$,$W_2$, $W_3$\} whose boundary
\{$\p W_1$, $\p W_2$, $\p W_3$\} is shown in Figure 9. Figure 9
also shows the longitude $\tilde \l$ of $\tilde M$, which is a
lift of $\l$.

This Heegaard splitting of $\tilde M$ is weakly reducible: $\p
X_4$ is disjoint from $\p W_3$. We now show that the closed, genus 2
surface $S$ obtained by compressing the Heegaard surface $\p \tilde H$ using
the disks $W_3$ and $X_4$ is essential
in $\tilde M$. It is enough to show that the surface $S$ is
incompressible in $\tilde M(2)$, which is the manifold obtained by
Dehn filling $\tilde M$ with the slope $2$. $\tilde M(2)$ has the
induced Heegaard splitting $\tilde H\cup \tilde C(2)$. Note that
$\tilde M(2)$ is the free 3-fold cyclic cover of $M(6)$, extending
the cover $\tilde M\ra M$, and that $\tilde C(2)$ is a handle body
of genus four covering the handle body $C(6)$ of genus two,
extending the cover $\tilde C\ra C$. Let $\tilde V$ be the filling
solid torus in $\tilde M(2)$ and let $W_4$ be a meridian disk of
$\tilde V$. Then $\{W_1, W_2, W_3, W_4\}$ is a disk system of the
handle body $\tilde C(2)$.

 Cutting $\tilde H$ along $X_4$, we get
a handle body $H_\#$ of genus three, and $\{X_1, X_2, X_3\}$ is a
disk system of $H_\#$. Using the Whitehead algorithm [S], we see
that $\p H_\#-\p W_3$ is incompressible in $H_\#$. In fact, from
Figure 9, we can read off the Whitehead graph of $\p W_3$ with
respect to the disk system $\{X_1, X_2, X_3\}$ of $H_\#$, which is
given as Figure 10. The graph is connected with no cut vertex,
which means, by the Whitehead algorithm, that $\p W_3$ must
intersect every essential disk of $H_\#$. Now by the Handle
Addition Lemma  due to Przytycki [P] and Jaco [J], the manifold
$H_\#\cup W_3\times I$, obtained by attaching the 2-handle
$W_3\times I$ to $H_\#$,  has incompressible boundary.

 \vspace{3mm}
\begin{figure}[!ht]
{\epsfxsize=5in \centerline{\epsfbox{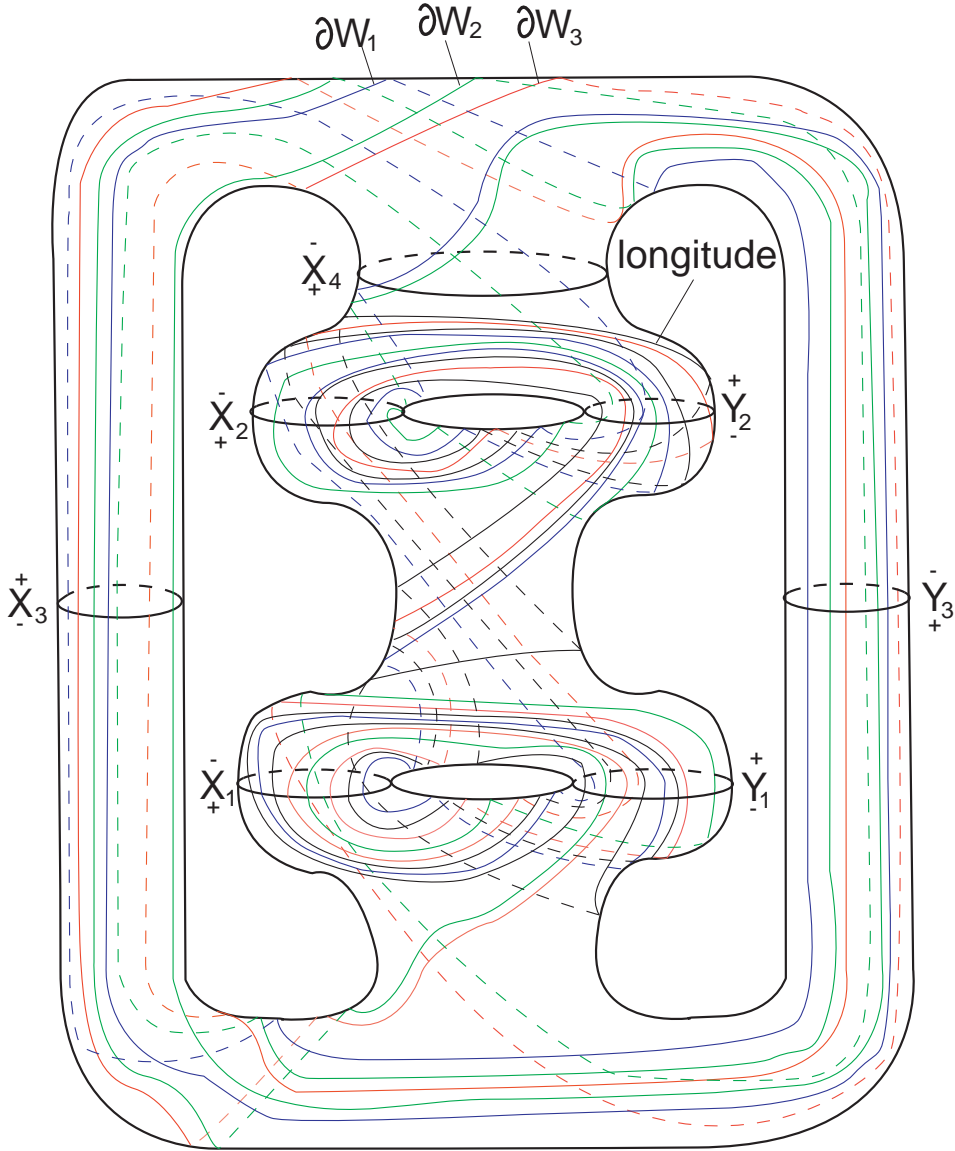}}\hspace{10mm}}
\caption{the Heegaard diagram of the 3-fold cyclic cover $\tilde
M$  and the longitude $\tilde \l$}
\end{figure}

 \vspace{3mm}
\begin{figure}[!ht]
{\epsfxsize=3in \centerline{\epsfbox{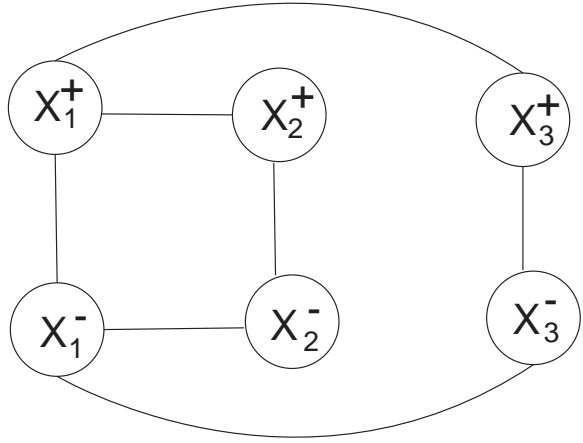}}\hspace{10mm}}
\caption{The Whitehead graph of $\p W_3$ with respect to the disk
system  $\{X_1,X_2,X_3\}$ of the handle body $H_\#$}
\end{figure}

On the other hand, cutting the handle body $\tilde C(2)$  along
the disk $W_3$,  we get a handle body $H_*$, which is homeomorphic
to $\tilde V$ with the two $1$-handles $W_1\times I$ and
$W_2\times I$ attached on
 $\p \tilde V$. The genus of $H_*$ is three, and
$\{W_1, W_2,W_4\}$ gives a disk system.
Let $\a\subset \p M$ be an essential simple closed curve of slope
$6$. We can easily see that with respect to the generators $x,y$
of $\pi_1(M)$, $$\a=\l x^6=yxyx^{-1}y^{-1}x^{-1}y^{-2}x^{-1}y^{-1}
x^{-1}yxyx^8.$$ Let $\tilde \a\subset \p\tilde M$ be a lift of $\a$.
Then $\tilde \a$ has slope $2$ in $\p\tilde M$ which can be
considered as the boundary of the disk $W_4$. Figure 11 shows
$\tilde \a=\p W_4$, $\p W_1$ and $\p W_2$ in $\p \tilde H$.

 \vspace{3mm}
\begin{figure}[!ht]
{\epsfxsize=5in \centerline{\epsfbox{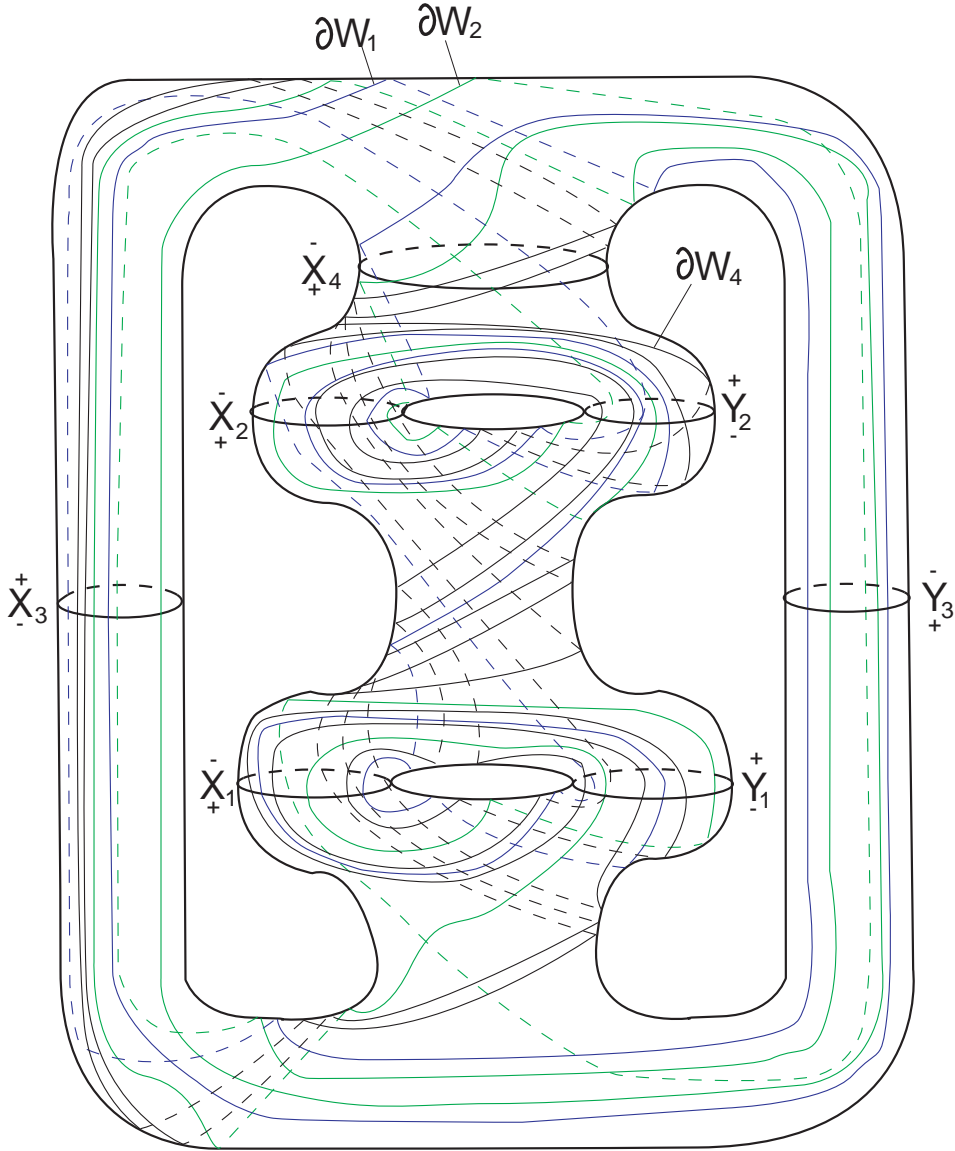}}\hspace{10mm}}
\caption{$\p W_4=\tilde \a$, $\p W_1$ and $\p W_2$ on the Heegaard
surface $\p \tilde H$}
\end{figure}

 \vspace{3mm}
\begin{figure}[!ht]
{\epsfxsize=3in \centerline{\epsfbox{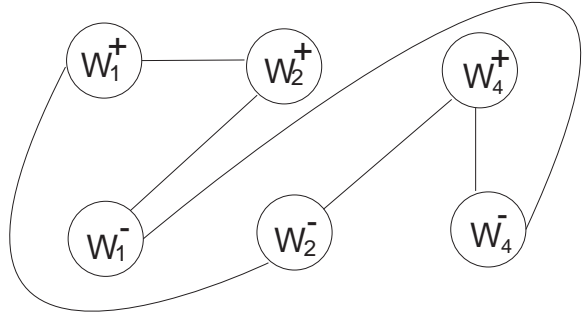}}\hspace{10mm}}
\caption{The Whitehead graph of $\p X_4$ with respect to the disk
system  $\{W_1,W_2,W_4\}$ of the handle body  $H_*$}
\end{figure}

  Again using the
Whitehead algorithm, we see that $\p H_*-\p X_4$ is incompressible
in $\tilde H_*$. In fact, from Figure 11, we can read off the
Whitehead graph of $\p X_4$ with respect to the disk system
$\{W_1, W_2, W_4\}$, which is given as Figure 12. The graph is
connected with no cut vertex, which means, by the Whitehead algorithm,
that $\p H_*-\p X_4$ is incompressible in $H_*$. Again by the
Handle Addition Lemma, the manifold $H_*\cup X_4\times I$ has
incompressible boundary of genus two. Note that $\p(H_* \cup
X_4\times I) =\p(\tilde H_\#\cup Y_3\times I)=S$ (up to a small
isotopy), and thus $S$ is incompressible in $\tilde M(2)$. But the
surface $S$ is contained $\tilde M$, and thus it is an
essential surface in $\tilde M$.

Obviously the longitude $\tilde \l$ in $\p \tilde M$ is isotopic
to an essential simple closed curve in the surface $S$, as shown
in Figure 9.
 The proof of Theorem \ref{3fold} is complete for $n=1$.

The proof for general $K_{2n+1}$, $n>0$, is similar.  The knot $K_{2n+1}$
is tunnel number one, with an unknotting tunnel shown in
Figure 2 (replacing the bottom three crossings by $2n+1$
crossings).  Let $M_n$ be the exterior of $K_{2n+1}$, $H'$ the
handlebody which is  a regular neighborhood of the knot and its
unknotting tunnel, $H=\overline{M_n-H'}$, and $D$ a meridian disk
of the unknotting tunnel. There is a corresponding
Heegaard splitting $M_n=H\cup_{\p H}C$, where $C$ is a compression
body. We let $\l$ be a standard longitude, and again we
deform the handlebody $H'$ by an isotopy in $S^3$ so
that its exterior $H$ can be recognized as a standard handlebody
in $S^3$, while tracing the corresponding
deformations of $\p D$ and $\l$ under the isotopy. We thus
get two essential disks $X$ and $Y$ which form a disk
system of $H$. From $\p D$, we get a geometric presentation of the
fundamental group $\pi_1(M_n)$ of $M_n$ with respect to the disk
system $\{X,Y\}$: $$\pi_1(M)=<x,y; (x^{-1}y^{-1})^{2n-1}
x^{-1}(yx)^{n+1} y^{-1}(x^{-1}y^{-1})^{2n-1}(xy)^{n+1}>. $$ Also
we get
$$\l=y(xy)^{n}(x^{-1}y^{-1})^nx^{-1}y^{-2}(x^{-1}y^{-1})^n
x^{-1}(yx)^{n+1}x.$$

 Let $\tilde M_n$ be the 3-fold cyclic cover of $M_n$
 and let $\tilde M_n=\tilde H\cup_{\p\tilde H}\tilde C$
 have the induced Heegaard splitting from that of
$M_n$, where $\tilde H$ is a genus four handle body which is the
corresponding 3-fold cyclic  cover of $H$ and $\tilde C$ a
compression body which covers $C$. Again the disk $X$ lifts to
three disks  $X_1, X_2, X_3$; and the disk $Y$ lifts to three
disks $Y_1, Y_2, Y_3$, as shown in Figure 9 (ignore the
$\p W_i$ and $\tilde \l$ part), and we pick the
meridian disk $X_4$ of $\tilde H$ as shown in Figure 9. Then
\{$X_1$, $X_2$, $X_3$, $X_4$\} forms a disk system of $\tilde H$.
The disk $D$ lifts to three disks \{$W_1$,$W_2$, $W_3$\} which
form a disk system of $\tilde C$. Again exactly one of the disks
\{$W_1$, $W_2$, $W_3$\}, say $W_3$, is disjoint from $X_4$, which
shows that  the Heegaard splitting of $\tilde M_n$  is weekly
reducible.  Again one can show that the surface $S$
obtained by compressing the Heegaard surface $\p \tilde H$ using
the disks $W_3$ and $X_4$ is an essential closed genus two surface
in $\tilde M_n$. In fact,
 cutting $\tilde H$ along $X_4$, we get
a handle body $H_\#$ of genus three and $\{X_1, X_2, X_3\}$ is a
disk system of $H_\#$. The  Whitehead graph of $\p W_3$ with
respect to the disk system $\{X_1, X_2, X_3\}$ of $H_\#$ is given
as Figure 13. The graph is connected with no cut vertex, which
means that $\p H_{\#}-\p W_3$ is incompressible. Thus by the
handle addition lemma, the manifold $H_\#\cup W_3\times I$,
obtained by attaching the 2-handle $W_3\times I$ to $H_\#$,  has
incompressible boundary.

 On the other hand, letting
$\tilde C(2)$ be the handle body obtained by Dehn filling $\tilde
C$ with slope $2$ and letting $W_4$ be a meridian disk of the filling
solid torus, then $\{W_1, W_2,W_3,W_4\}$ forms a disk system of
$\tilde C(2)$. Cutting  $\tilde C(2)$ along the disk $W_3$, we get
a handlebody $H_*$ with disk system $\{W_1,W_2,W_4\}$. Let
$\a\subset \p M$ be an essential simple closed curve of slope $6$.
Then with respect to the generators $x,y$ of $\pi_1(M)$, $$\a=\l
x^6=y(xy)^{n}(x^{-1}y^{-1})^nx^{-1}y^{-2}(x^{-1}y^{-1})^n
x^{-1}(yx)^{n+1}x^6.$$ We may consider $\p W_4$ as a lift of $\a$.
From the word $\a$, we can draw $\p W_4$ on $\p \tilde H$.
Consequently we can read off the Whitehead graph of $\p X_4$ with
respect to the disk system $\{W_1, W_2, W_4\}$ and see that the
graph is the same as that shown in Figure 12, showing that  $\p
H_*-\p X_4$ is incompressible in $H_*$. Thus the manifold
$H_*\cup X_4\times I$ has incompressible boundary of genus two. We
thus have justified the incompressibility of the surface $S$ in
$\tilde M_n(2)$ and thus in $\tilde M_n$.

 \vspace{3mm}
\begin{figure}[!ht]
{\epsfxsize=5in \centerline{\epsfbox{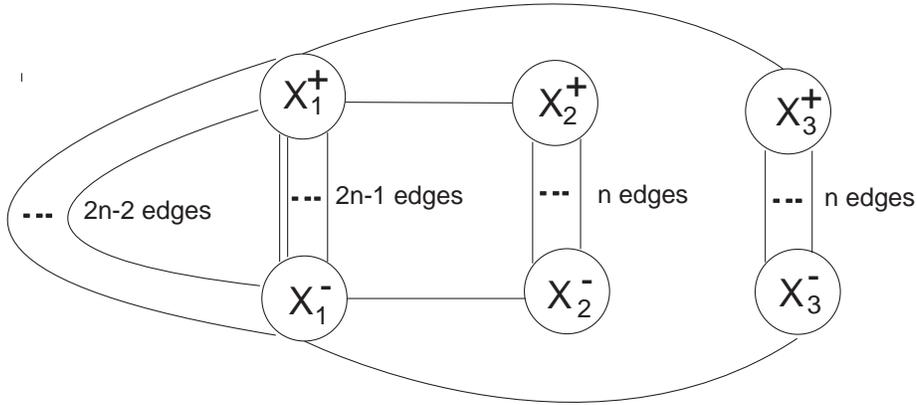}}\hspace{10mm}}
\caption{The Whitehead graph of $\p W_3$ with respect to the disk
system  $\{X_1,X_2,X_3\}$ of the handle body  $H_\#$}
\end{figure}

Finally  the longitude $\tilde \l$ in $\p \tilde M$ is isotopic to
an essential simple closed curve in the surface $S$, which is
obvious.
 The proof for the general case is
complete. \qed

\vspace{3mm}
\begin{figure}[!ht]
{\epsfxsize=3in \centerline{\epsfbox{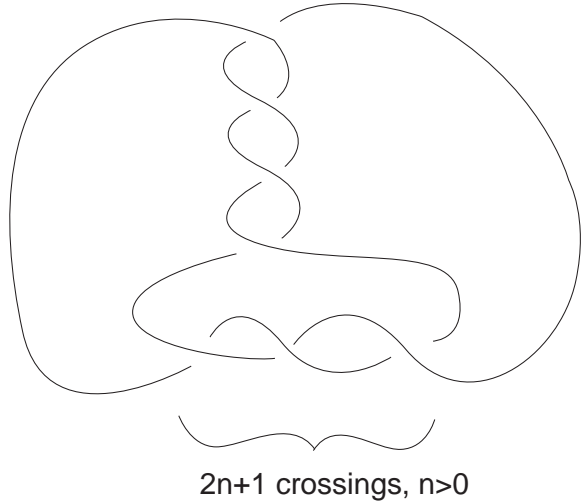}}\hspace{10mm}}
\caption{The  knot $J_{2n+1}$}
\end{figure}

Let $J_{2n+1}$, $n>0$,  be the family of two bridge knots shown in
Figure 14. Note that these knots are hyperbolic, small and
non-fibered with genus two Seifert surfaces.

\begin{thm}\label{j}
 The $5$-fold cyclic cover of the exterior of $J_{2n+1}$
is large  and every Dehn filling of the exterior of
$J_{2n+1}$  with slope
$5p/q$, $(5p,q)=1$, $|p|>1$, yields  a virtually Haken 3-manifold,
for every $n>0$.
\end{thm}

This theorem gives another family of non-Haken, virtually Haken
 3-manifolds to which the results of [CL] do not apply.
 As the proof of Theorem \ref{j} is very similar to that of
 Theorem \ref{3fold}, we omit the details and indicate only the steps.
 In fact the exterior of $J_{2n+1}$ is tunnel number one
 and a genus two Heegard splitting of it can be explicitly
 given as
 in the case for the exterior of the twist knot $K_{2n+1}$.
 In the $5$-fold cyclic cover of the exterior of $J_{2n+1}$,
 the lifted Heegaard surface
 is of genus $6$ and  can be compressed   along two reducing disks,
 one on each side of the Heegaard surface,
 to a closed incompressible surface of genus $4$.
 Also a lift of the longitude can be isotoped into the
 resulting incompressible surface.

 \vspace{3mm}
\begin{figure}[!ht]
{\epsfxsize=4.5in \centerline{\epsfbox{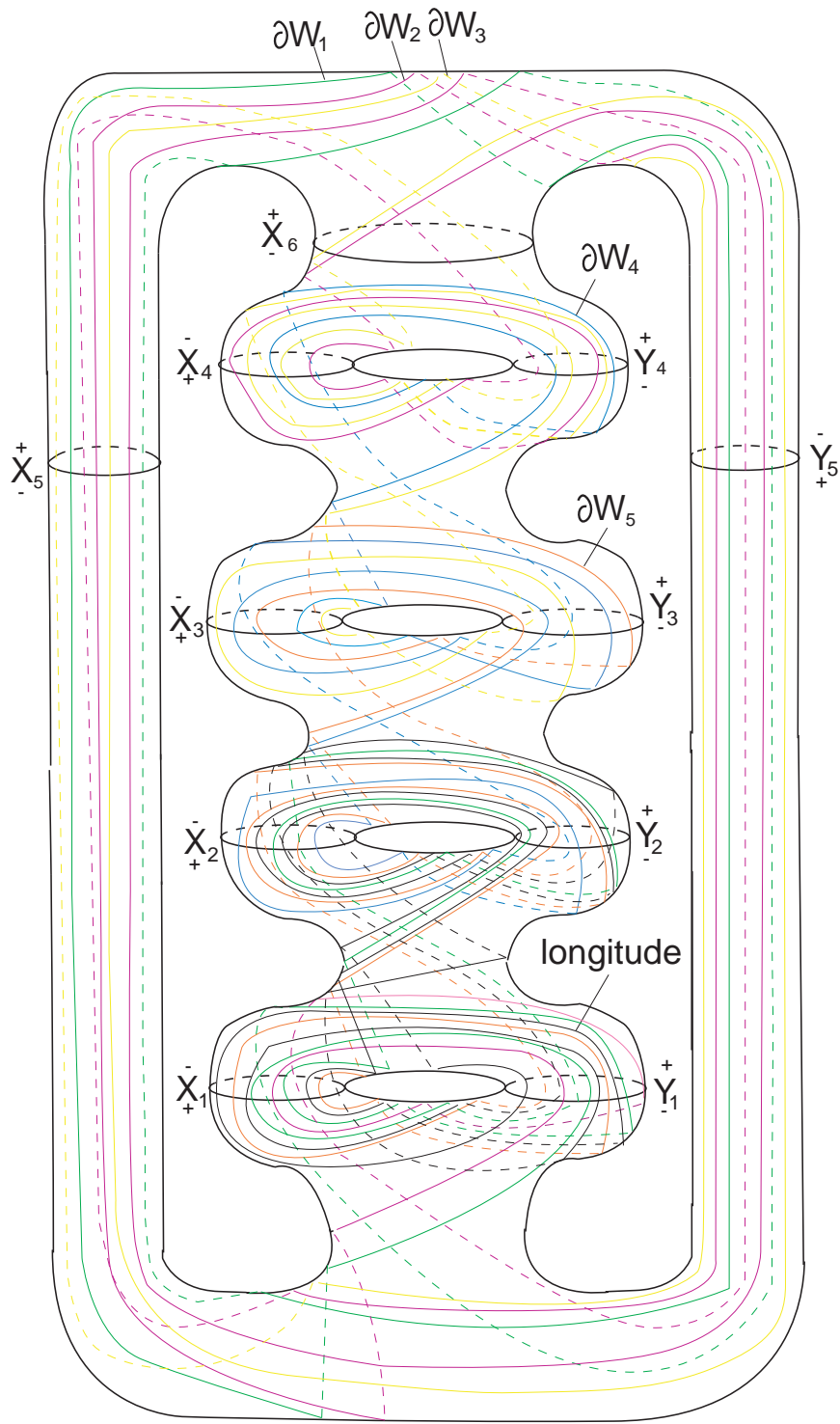}}\hspace{10mm}}
\caption{The Heegaard splitting of the 5-fold cover of $M$}
\end{figure}

We now go back to the twist knots $K_{2n+1}$
and prove the following Theorem \ref{5fold}.
Although the result of the theorem is covered by [CL],
we have included it primarily because its proof illustrates
 two complications which arise in more general settings.
 First, we have to deal with multi
 2-handle additions, which requires the multi 2-handle
 addition theorem of Lei [L]. Also,
 one of the Whitehead graphs contains a cut vertex,
 and must be simplified using Whitehead moves.

 \vspace{3mm}
\begin{figure}[!ht]
{\epsfxsize=3in \centerline{\epsfbox{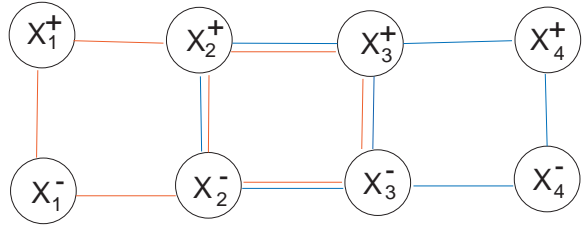}}\hspace{10mm}}
\caption{The Whitehead graph of \{$\p W_4$,$\p W_5$\} with respect
to the disk system $\{X_1,X_2,X_3,X_4\}$ of the handle body
$H_\#$}
\end{figure}

\begin{thm}\label{5fold}
 The $5$-fold cyclic cover of the exterior $M_n$ of $K_{2n+1}$
is large for every  $n>0$. Every Dehn filling of $M_n$  with slope
$5p/q$, $(5p,q)=1$, $|p|>1$, yields  a virtually Haken 3-manifold.
\end{thm}

\pf  Again we give details only for the $n=1$ case.
 We continue to use the
Heegaard splitting of $M=M_1=H\cup C$ as given in the proof of Theorem
\ref{3fold}.  Let $\tilde M$ be the 5-fold cyclic cover of $M$
with the induced Heegaard splitting from that of $M$. The Heegaard
diagram of $\tilde M$ is shown in Figure 15. The genus six handle
body of Figure 14 is $\tilde H$ which covers $H$. The disks $X$
and $Y$ of $H$ lift to  disks $X_1,..., X_5$ and  $Y_1, ..., Y_5$,
as shown in Figure 15. Pick the meridian disk $X_6$ of $\tilde H$
as shown in Figure 15. Then \{$X_1 ,X_2,X_3,X_4,Y_5, X_6$\} forms
a disk system of $\tilde H$. The disk $D$ lifts to five disks
\{$W_1,W_2,W_3, W_4,W_5$\} whose boundaries  are shown in Figure 15. Figure 15
also shows a longitude $\tilde \l$ of $\tilde M$, which is a lift
of the longitude $\l$ of $M$.

 \vspace{3mm}
\begin{figure}[!ht]
{\epsfxsize=4.5in \centerline{\epsfbox{f17.ai}}\hspace{10mm}}
\caption{$\p W_6=\tilde \a$, $\p W_1$, $\p W_2$, $\p W_3$ on the
Heegaard surface $\p \tilde H$}
\end{figure}

This Heegaard splitting of $\tilde M$ is weakly reducible: \{$\p
Y_5$, $\p X_6$\} is disjoint from \{$\p W_4$, $\p W_5$\}. We now
show that the surface $S$ obtained by compressing the Heegaard
surface $\p \tilde H$ using these four disks is an essential
closed genus two surface in $\tilde M$. It is enough to show that
the surface $S$ is incompressible in $\tilde M(2)$, which is the
free 5-fold cyclic cover of $M(10)=H\cup C(10)$, and has the
induced Heegaard splitting $\tilde H\cup \tilde C(2)$. Let $\tilde V$ be
the filling solid torus in $\tilde M(2)$ and let $W_6$ be a meridian
disk of $\tilde V$. Then $\{W_1, ..., W_5, W_6\}$ is a disk system
of the handle body $\tilde C(2)$.

 \vspace{3mm}
\begin{figure}[!ht]
{\epsfxsize=3in \centerline{\epsfbox{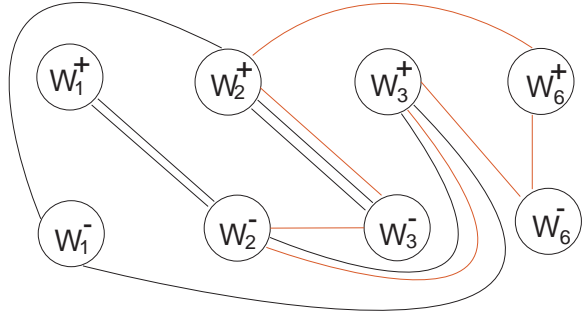}}\hspace{10mm}}
\caption{The Whitehead graph of \{$\p Y_5$, $\p X_6$\} with
respect to the disk system  $\{W_1,W_2,W_3,W_6\}$ of the handle
body $H_*$}
\end{figure}

 Cutting $\tilde H$ along $Y_5, X_6$, we get
a handle body $H_\#$ of genus four and $\{X_1, X_2, X_3, X_4\}$ is
a disk system of $H_\#$. The Whitehead graph of \{$\p W_4$, $\p
W_5$\} with respect to the disk system $\{X_1, ..., X_4\}$ of
$H_\#$ is given in Figure 16. The graph is connected with no cut
vertex, which means that the surface $\p H_\#-\{\p W_4, \p W_5\}$
is incompressible in $H_\#$. Moreover as $\p W_4$ is disjoint from
the disk $X_1$, and $\p W_5$ is disjoint from the disk $X_4$, each
of  the surfaces $\p H_\#-\p W_4$ and $\p H_\#-\p W_5$ is
compressible in $H_\#$. Therefore all the conditions of the
multi-handle addition theorem of [L] are satisfied, and thus the
manifold $H_\#\cup W_4\times I\cup W_5\times I$ has incompressible
boundary.

 \vspace{3mm}
\begin{figure}[!ht]
{\epsfxsize=3in \centerline{\epsfbox{f19.ai}}\hspace{10mm}}
\caption{(a) The resulting graph after the Whitehead move with
respect to the cut vertex $W_2^-$ of Figure 18. (b) The resulting
graph after the Whitehead move with respect to the cut vertex
$W_3^-$ of part (a) }
\end{figure}

On the other hand, cutting the handle body $\tilde C(2)$  along
the disks $W_4$ and $W_5$,  we get a handle body $H_*$, with disk system
$\{W_1,W_2, W_3, W_6\}$. Let $\a\subset
\p M$ be an essential simple closed curve of slope $10$. Then
$$\a=\l x^{10}=yxyx^{-1}y^{-1}x^{-1}y^{-2}x^{-1}y^{-1}
x^{-1}yxyx^{12}.$$ Let $\tilde \a\subset \p\tilde M$ be a lift
$\a$. Then $\tilde \a$, which can be
considered as the boundary of the disk $W_6$,
 has slope $2$ in $\p\tilde M$. Figure 17 shows
$\tilde \a=\p W_6$, $\p W_1$,$\p W_2$, $\p W_3$ in $\p \tilde H$.

From Figure 17, we can read  off the Whitehead graph of \{$\p
Y_5$, $\p X_6$\} with respect to the disk system $\{W_1, W_2, W_3,
W_6\}$ of $H_*$, which is given as Figure 18. The graph is
connected but has a cut vertex (the vertex $W_2^-$). Applying
Whitehead moves to the graph twice with results shown in Figure
19, we end up with a graph (shown in Figure 19 (b)) which is
connected with no cut vertex. This means that the surface $\p
H_*-\{\p Y_5\cup \p X_6\}$ is incompressible in $H_*$. From Figure
16, we also see that
 $\p Y_5$ is disjoint from $\p W_6$ and $\p X_6$ is
disjoint from $\p W_1$.
Thus each
of  the surfaces $\p H_*-\p Y_5$ and $\p H_*-\p X_6$ is
compressible in $H_*$. Again  the
multi-handle addition theorem of [L] implies
that  the manifold $H_*\cup X_6\times I\cup Y_5\times I$
 has incompressible boundary.
Therefore the genus two surface
 $S=\p(H_*\cup X_6\times I\cup Y_5\times I)=
\p (H_\#\cup W_4\times I\cup W_5\times I)$
is  incompressible in $\tilde M(2)$ and thus
is essential in $\tilde M$.

Obviously $\tilde \l$ can be isotoped into $S$. The proof of
Theorem \ref{5fold} is complete in case $n=1$. The proof for the
general case is similar (cf the proof of Theorem \ref{3fold}
in general case). We leave the details to the  reader to verify.
\qed

\newpage 
\begin{center}
{\bf References}
\end{center}

\noindent [BZ]  S. Boyer and X. Zhang,
 Virtual Haken 3-manifolds and Dehn filling,
Topology 39 (2000) 103-114.

\noindent
[CG] A. Casson and C. Gordon,
 Reducing Heegaard splittings. Topology Appl. 27 (1987) 275--283.

\noindent
[CGLS] M. Culler, C. Gordon, J. Luecke and P. Shalen,
Dehn surgery on knots, Ann. of Math. 125 (1987) 237-300.

\noindent [CL] D. Cooper and D. Long,  Virtually Haken
Dehn-filling, J. Differential Geom. 52 (1999) 173--187.

\noindent [DT] N. Dunfield and W. Thurston, The virtusl Haken
conjecture: experiments and examples, preprint.

\noindent
[HT] A. Hatcher and W. Thurston, Incompressible surfaces
in $2$-bridge knot complements. Invent. Math. 79 (1985) 225--246.

\noindent
[J] W. Jaco, Adding a $2$-handle to a $3$-manifold: an
application to property $R$. Proc. Amer. Math. Soc. 92 (1984)
288--292.

\noindent [K] R. Kirby, Problems in Low-Dimensional Topology,
Geometric topology. Edited by William H. Kazez. AMS/IP Studies in
Advanced Mathematics, 2.2.  (1997) 35-473.

\noindent
[La] M. Lackenby,
 Heegaard splittings, the virtually Haken conjecture and Property tau,
 preprint.

\noindent
 [L] F. Lei, A proof of Przytycki's conjecture on
n-relator 3-manifolds, Topology (1995) 473-476.

\noindent [M] J. Masters, Virtual homology of surgered torus
bundles. Pacific J. Math. 195 (2000) 205--223.

\noindent
[P] J.  Przytycki,
 Incompressibility of surfaces after Dehn surgery.
  Michigan Math. J. 30 (1983) 289--308.

\noindent
[S] J. Stallings, Whitehead graphs on handlebodies,
preprint.

\vspace{10mm} \noindent Mathematics Department, Rice University,
Houston, TX 77005\newline  mastersj@rice.edu

\vspace{5mm} \noindent Mathematics Department, SUNY at Buffalo,
Buffalo, NY 14260-2900\newline menasco@math.buffalo.edu

\vspace{5mm} \noindent Mathematics Department, SUNY at Buffalo,
Buffalo, NY 14260-2900\newline xinzhang@math.buffalo.edu
\end{document}